\documentclass[12pt,oneside]{amsart}
\usepackage{amscd,amsmath,amssymb,amsfonts,a4}
\newlength{\rulebreite}

\def\timesover#1#2#3{\ \xymatrix@1@=0pt@M=0pt{ _{#1}&\times&_{#2} \\& ^{#3}&}\ }
\def\otimesover#1#2#3{\ \xymatrix@1@=0pt@M=0pt{ _{#1}&\otimes&_{#2} \\& ^{#3}&}\ }
\usepackage[all]{xy}
\theoremstyle{plain}
\newtheorem{thm}{Theorem}
\newtheorem{lem}[thm]{Lemma}
\newtheorem{cor}[thm]{Corollary}
\newtheorem{prop}[thm]{Proposition}
\newtheorem{prop-ex}[thm]{Example}
\theoremstyle{definition}

\newtheorem{claim}[thm]{Claim}
\numberwithin{thm}{section}
\numberwithin{equation}{section}
\newcommand{\eq}[2]{\begin{equation}\label{#1}#2 \end{equation}}
\newcommand{\ml}[2]{\begin{multline}\label{#1}#2 \end{multline}}
\newcommand{\ga}[2]{\begin{gather}\label{#1}#2 \end{gather}}
\newcommand{\inj}{\hookrightarrow}

\newcommand{\Pic}{{\rm Pic}}

\newcommand{\Spec}{{\rm Spec \,}}

\newcommand{\sL}{{\mathcal L}}

\newcommand{\sO}{{\mathcal O}}

\newcommand{\A}{{\mathbb A}}

\newcommand{\C}{{\mathbb C}}

\newcommand{\G}{{\mathbb G}}

\newcommand{\BL}{{\mathbb L}}

\newcommand{\BP}{{\mathbb P}}
\newcommand{\Q}{{\mathbb Q}}

\newcommand{\Z}{{\mathbb Z}}

\begin{document}
\title[Rationality]{On a rationality question in the Grothendieck ring of varieties}
\author{H\'el\`ene Esnault}
\address{
Universit\"at Duisburg-Essen, Mathematik, 45117 Essen, Germany}
\email{esnault@uni-due.de}
\author{Eckart Viehweg}
\address{Universit\"at Duisburg-Essen, Mathematik, 45117 Essen, Germany}
\email{viehweg@uni-due.de}
\date{October 25, 2009}
\thanks{Partially supported by  the DFG Leibniz Preis, the SFB/TR45, the ERC Advanced Grant 226257}
\begin{abstract} We discuss elementary rationality questions in the Grothendieck ring of varieties for the quotient of a finite dimensional vector space over a characteristic 0 field by a finite group.
\end{abstract}
\maketitle
\section{Introduction} 
Let $k$ be a field. One defines the {\it Grothendieck group of varieties} $K_0({\rm Var}_k)$ {\it  over } $k$ \cite[Definition~2.1]{Ni} to be the free abelian group generated by $k$-schemes modulo the subgroup spanned the scissor relations 
$$[X]=[X\setminus Z]+ [Z]$$ 
where $Z\subset X$ is a closed subscheme. The product 
$$[X\times_k Y]=[X]\cdot[Y]$$
for two $k$-schemes makes it a commutative ring, with unit $1=[\Spec k]$. As the underlying  topological space of the complement $X\setminus X_{{\rm red}}$ is empty, $[X]=[X_{{\rm red}}]$. This justifies the terminology  ``varieties'' rather than ``schemes''.\\[.1cm]
In characteristic 0, first examples of 0-divisors in this ring were shown to exist by Poonen \cite{Poo}. He constructed two abelian varieties $A, B$ over $\Q$ such that 
$$0=([A]-[B])\cdot ([A]+[B]) \in K_0({\rm Var}_{\Q})$$
but with $$[A\otimes_{\Q} k]\neq [B\otimes_{\Q} k] \in K_0({\rm Var}_k)$$ for all field extensions $\Q\hookrightarrow{} k$. 
The main tool to distinguish those two classes relies ultimately on a deep insight in the structure of birational morphisms, gathered in the {\it Weak Factorization Theorem} \cite{AKMW}. It implies both the presentation of $K_0({\rm Var}_k)$ as the free group generated by smooth projective varieties modulo the blow up relation \cite{Bi} and the isomorphism $K_0({\rm Var}_k)/\langle \BL \rangle \xrightarrow{\cong} \Z[SB]$ \cite{LaLu}.  Here $\BL$ is the class of the affine line $\A^1$ over $k$, $\langle \BL \rangle$ is the ideal spanned by it, $\Z[SB]$ is the free abelian group on stably birational classes of projective smooth $k$-varieties, endowed with the ring structure stemming from the product of varieties over $k$. So there are no relations in $\Z[SB]$ and this allows to recognize certain classes. Of course this does not help in understanding $\BL$, and the question whether or not $\BL$ is a 0-divisor remains open. \\[.1cm]
Later Koll\'ar \cite{Ko} used $\Z[SB]$ to distinguish in characteristic 0 the $K_0({\rm Var}_k)$-classes of non-trivial Severi-Brauer varieties 
from trivial ones. R\"okaeus \cite{Ro} and Nicaise \cite{Ni}, using in addition specialization of $K_0({\rm Var}_k)$ from $k$ to  finite fields, studied 0-divisors which are classes of 0-dimensional varieties, in particular those of the form $\Spec K$ for a non-trivial field extension of a number field $k$. This indicates that one can not
expect ``descent''. For two $k$-varieties $X$ and $Y$ the equality
$$[X\times_k \Spec K] = [Y\times_k \Spec K] \in K_0({\rm Var}_K)$$
implies, by the projection formula, that
$$[X]\cdot [\Spec K] = [Y]\cdot [\Spec K] \in K_0({\rm Var}_k).$$
Indeed, the field extension $\iota: k\hookrightarrow K$, induces a ring homomorphism 
$$\iota^*:
K_0({\rm Var}_k)\to K_0({\rm Var}_K),$$ defined by $\iota^*[X]=[X\times_k \Spec K]$, while $$\iota_*:   K_0({\rm Var}_K) \to K_0({\rm Var}_k)$$ is the homomorphism of abelian groups which takes the class of the $K$-variety $[Z]$
to the same class viewed as a $k$-variety. The projection formula says $$\iota_*\big(\iota^*[X]\cdot [Z]\big)=[X]\cdot \iota_*[Z],$$ 
thus applied to $[Z]=1=[\Spec K]\in K_0({\rm Var}_K)$ it yields the formula $$\iota_*[X\times_k K]=[X]\cdot [\Spec K] \in K_0({\rm Var}_k).$$  
However, the relation $[X]\cdot [\Spec K] = [Y]\cdot [\Spec K] \in K_0({\rm Var}_k)$
does not imply the equality $[X]=[Y]\in K_0({\rm Var}_k)$.\\[.1cm]
For applications of the Grothendieck ring, it is of importance to understand the class of quotients $[X/G]$ where $X$ is a variety and $G$ is a finite group acting on it. 
In \cite[Lemma~5.1]{Loo}, Looijenga shows that if $k$ is an algebraically closed field of characteristic 0, and if $G$ is a finite abelian group acting linearly on a finite dimensional $k$-vector space $V$, then 
\ga{eq1.1}{[V/G]=\BL^{{\rm dim}_k V} \in K_0({\rm Var}_k).}
In fact the formula \eqref{eq1.1}, as well as its proof, remain valid if
$k$ is any field of characteristic 0 containing the $|G|$-th roots of $1$.
However the condition that $G$ be abelian is  essential, as shown by Ekedahl. Indeed, \cite[Proposition~3.1,~ii)]{Ek} together with \cite[Corollary~5.2]{Ek} show that for $G \subset GL(V), \ V\cong \C^n$ as in Saltman's example \cite{S}, the
class of   $\lim_{m\to \infty} [V^m/G]/\BL^{nm}$ in the completion $ \widehat{ K_0({\rm Var}_{\C})}$ of $K_0({\rm Var}_{\C})[\BL^{-1}]$ by the dimension filtration, is not equal to $1$. This implies in particular  that for $m$ large enough, 
  $\BL^{nm}\neq [V^m/G]\in K_0({\rm Var}_{\C})$.
   \\[.2cm] 
In this note, we discuss possible simple generalizations of Looijenga's formula in various ways. 
Our first result is the following.
\begin{lem}  \label{thm1.1} Let $G$ be a finite abelian group with quotient
$G \to \Gamma$. Let $k$ be a field of characteristic 0 and let $K\supset k$ be an abelian Galois extension with Galois group $\Gamma$. 
Assume, that the Galois action of $\Gamma$ on $K$ lifts to a $k$-linear action of
$G$ on a finite dimensional $K$-vector space $V$. If, for $N={\rm exp}(G)$, all $N$-th roots of $1$ lie in $k$, then \eqref{eq1.1} holds, i.e.
$$[V/G]= {\mathbb L}^{{\rm dim}_K V} \in K_0({\rm Var}_k).$$
\end{lem}
The condition that $k$ contains the $N$-th roots of $1$ is really necessary.
In particular, if one allows the group $G$ to act non-trivially on the ground field,
the equation \eqref{eq1.1} is not compatible with descent to smaller ground fields.  
\begin{prop-ex} \label{ex1.2}
Assume $k=\Q, \ K=\Q(\sqrt{ -1})$, $V=K\otimes_{\Q}\Q^2 $, and let $G$ be the subgroup of the group of $\Q$-linear automorphisms of $V$ spanned by  
\ga{}{\sigma= \begin{pmatrix} 1 & 0\\
0 & -1 \end{pmatrix}\otimes\begin{pmatrix} 0 & 1\\
-1 & 0\end{pmatrix}   \notag }
where the chosen basis of $K$ as a 2-dimensional vector space over $\Q$ is $(1, \sqrt{ -1})$. The group  $G$ is cyclic of  order $4$ and
$$\BL^2\neq [V/G] \in K_0({\rm Var}_\Q).$$
\end{prop-ex}
If $G\subset GL_k(V)$ is a finite group acting linearly on a finite dimensional vector space $V$ over a characteristic 0 field $k$, then $G$ acts semi-simply. So as a $G$-representation,  $V=\bigoplus_i V_i\otimes T_i$, where $V_i$ is an irreducible representation with ${\rm Hom}_G(V_i, V_j)=\delta_{ij}\cdot k$, and $T_i$ is the trivial representation of dimension $m_i$ equal to the multiplicity of $V_i$ in $V$.
If $G$ is commutative of exponent $N$ and if the $N$-th roots of $1$ lie in $k$, then $d_i={\rm dim}_k V_i=1$.
Since $V_i/G$ is normal and one dimensional, it is smooth. So the starting point of Looijenga's proof of \eqref{eq1.1} is the simple observation that there is a $k$-isomorphism $V_i/G\cong V_i$ of $k$-varieties. The proof of \eqref{eq1.1} then proceeds by stratifying $V$.\\[.1cm]
For $d_i\ge 2$, the quotient $V_i/G$ might be singular, thus it can not be isomorphic
to $V_i$, not even over a field extension. Nevertheless, one can show that 
the formula \eqref{eq1.1} remains true for irreducible two dimensional representations, or
after stratifying, whenever all the $d_i$ are 1 or 2 and $G$ is a prime power order cyclic group.
\begin{prop} \label{prop1.3}
Let $k$ be a field of characteristic 0 and let $V$ be a finite dimensional $k$-vector space. Let $G\to GL_k(V)$ be a linear representation of a finite abelian group.
\begin{itemize}
\item[1)] If $\dim_k V\le 2$, then \eqref{eq1.1} holds true. 
\item[2)] If $G$ is cyclic of prime power order, and if each irreducible
subrepresentation $V_i$  has $\dim(V_i) \le 2$,  then \eqref{eq1.1} holds true. 
\end{itemize}
\end{prop}
The main reason for the restriction to $\dim(V_i) \le 2$ is that in this case $\BP(V_i)\cong \BP^1_k$ and hence $\BP(V_i)/G\cong \BP^1_k$ as well. If $V$ is an irreducible representation of dimension $d \ge 3$ a similar statement fails, and we were unable to prove the equation \eqref{eq1.1}.\\[.2cm]
{\it Acknowledgments:} We thank Johannes Nicaise for several discussions which motivated this little elementary note. Lemma \ref{thm1.1} and Example \ref{ex1.2} were shown in a letter to him dated September 27, 2008. They have been further worked out by Tran Nguyen Khanh Linh and Le Hoang Phuoc in their master thesis in Essen, 2009. In addition, we thank the referee for pointing out the reference \cite{Ek} to us and for a careful reading which helped us to improve the exposition of the manuscript.

\section{Proof of Lemma \ref{thm1.1}}
By assumption $G\subset GL_k(V)$ lifts the action of the quotient $\Gamma$ on $K$,
hence writing 
\ga{}{1 \longrightarrow{} H \longrightarrow{} G \xrightarrow{ \ \varrho \ } \Gamma \longrightarrow{} 1,\notag}
one has $\sigma(\lambda\cdot v)=\gamma(\lambda)\cdot \sigma(v)$, for all $\sigma\in G, \ \gamma=\varrho(\sigma)$, for all $\lambda \in K$ and for all $v\in V$. 
In particular $H$ is a subgroup of $GL_K(V)$. This defines the fiber square 
\ga{eq2.1}{ \xymatrix{ \ar @{} [dr] |{\square}\ar[d] V/H \ar[r] & \Spec K \ar[d]\\
V/G \ar[r] & \Spec k.}
}
By the rationality assumption, $\mu_N(k)\cong_k \Z/N$, for $N={\rm exp}(G)$, and hence 
the characters of $G$ are $k$-rational. So writing $\hat{H}$ for the character group of $H$
and $V_\chi(H)$ for the eigenspace with respect to the character $\chi$ of $H$,
one has a fortiori the $K$-eigenspace decomposition 
\ga{}{   V=\bigoplus_{\chi \in \hat{H}} V_{\chi}(H).\notag}
Since $G$ is commutative the subspace $V_{\chi}(H)$ of $V$ is $G$-invariant.\\[.1cm]
Now on the geometric side, one proceeds as in Looijenga's Bourbaki lecture \cite[Lemma~5.1]{Loo}. Write 
\ga{}{V=\prod_{\chi \in \hat{H}} V_{\chi}(H) \notag}
for the product as $K$-schemes. For $\{0\}=\Spec K$ one sets $V^\times_{\chi} = V_{\chi}(H) \setminus \{0\}$ and defines the stratification
\ga{eq2.2}{V=\bigsqcup_{I \subset \hat{H}} V_I, \mbox{ \ \ with \ \ } V_I=\prod_{\chi \in I} V^\times_{\chi}. }
The product in \eqref{eq2.2} is defined over $K$. The $\G_m$-fibration $V^\times_{\chi} \to \BP(V_{\chi}(H))$ is the
structure map of the geometric line bundle 
$\sO_{\BP(V_{\chi}(H))}(-1)$, restricted to the complement of the zero-section. 
It is defined over $K$ and $G$-equivariant. The subgroup $H$ acts trivially on
$\BP(V_{\chi}(H))$ and by multiplication with $\chi$ on the geometric fibres
of $V^\times_{\chi} \to \BP(V_{\chi}(H))$.\\[.1cm]
So for $I \subset \hat{H}$ given, the $K$-morphism 
\ga{}{ V_I\to \prod_{\chi \in I} \BP(V_{\chi}(H)) \notag}
is a $G$-equivariant fibration, locally trivial for the Zariski topology. The fibres 
are isomorphic to $\G_m^{\#I}\cong \prod_{\chi\in I} \G_{m,\chi}$, with $\G_{m,\chi}\cong \G_m$, hence 
\eq{eq2.3}{[V_I]=[\G_m^{\#I}]\cdot \prod_{\chi\in I} [ \BP^{r_\chi}_K]\mbox{ \ \ in \ \ }K_0({\rm Var}_K).}
The action of $H$ is trivial on $\prod_{\chi\in I} \BP(V_{\chi}(H))$ and
on the factor $\G_{m,\chi}$ of $\G_m^{\#I}$ the group $H$ acts by multiplication
with $\chi$. One obtains an induced $K$-morphism 
\ga{}{ V_I/H\to \prod_{\chi \in I} \BP(V_{\chi}(H)) \notag}
which is still a Zariski locally trivial fibration with fibre
\eq{eq2.4}{\G_m^{\#I}\cong\big[\prod_{\chi\in I} \G_{m,\chi}\big]/H.}
Since $G$ respects the
decomposition $V_I=\prod_{\chi\in I} V_{\chi}(H)$ and since the $G$ action on
$\BP(V_{\chi}(H))$ factors through $\Gamma$ one
finds
$$
\big[ \prod_{\chi \in I} \BP(V_{\chi}(H))\big]/G = \big[ \prod_{\chi \in I} \BP(V_{\chi}(H))\big]/\Gamma = \prod_{\chi \in I} \big(\BP(V_{\chi}(H))/\Gamma\big),
$$
where the first two products are over $K$ whereas the one on the right hand side is the product of $k$-varieties. 
Remark that $\BP(V_{\chi}(H))/G=\BP(V_{\chi}(H))/\Gamma$ is a $k$-form of $\BP^{r_\chi}_k$ for
$r_\chi={\rm dim}_K V_{\chi}(H) -1$.  \\[.1cm]
The fiber square \eqref{eq2.1} is the composite of two fibre squares
\ga{eq2.5}{ \xymatrix{ \ar @{} [dr] |{\square}\ar[d] V_I/H \ar[r] & \ar @{} [dr] |{\square}
\ar[d] \prod_{\chi \in I} \BP(V_{\chi}(H)) \ar[r] &
 \Spec K \ar[d]\\
V_I/G \ar[r] & \prod_{\chi \in I} \big(\BP(V_{\chi}(H))/\Gamma\big) \ar[r] & \Spec k.}
}
\begin{claim} \label{claim2.1} The $k$-form $\BP(V_{\chi}(H))/\Gamma$ of $\BP^{r_\chi}_k$ is split, the $k$-morphism 
$$
V_I/G\to \prod_{\chi \in I} \BP(V_{\chi}(H))/\Gamma$$ 
is a $\G_m^{\#I}$-fibration, locally trivial for the Zariski topology, and 
hence 
$$
[V_I/G]=[\G_m^{\#I}]\cdot \prod_{\chi\in I} [ \BP^{r_\chi}_k]\mbox{ \ \ in \ \ }K_0({\rm Var}_k).
$$ 
\end{claim}
\begin{proof}
By assumption $k$ contains the $N$-th roots of 1 for $N={\rm exp}(G)$ and hence the
characters $\chi\in \hat{H}$ are defined over $k$. \\[.1cm]
Then $V_{\chi}(H)$, regarded as a $k$-vector space, has a $G$-eigenvector $v$.
The line $\langle v \rangle_K$ defines a point $c\in \BP(V_\chi(H))(K)$. Since the action of $G$ on $K(c)=K$ factors through the Galois action of $\Gamma$ on $K(c)$, the image of $c$ lies in $(\BP(V_\chi(H))/G)(k)$.\\[.1cm]
In addition, in \eqref{eq2.4} the action of $H$ on
$\prod_{\chi\in I} \G_{m,\chi}$ is given by multiplication with
$\chi$, hence it is defined over $k$. Then $\big[\prod_{\chi\in I} \G_{m,\chi}\big]/H$
is obtained by base extension from a $k$-variety, isomorphic to $\G_{m}^{\#I}$.\\[.1cm]
Using that the left hand side of \eqref{eq2.5} is a fibre product and that
$$
V_I/H \to \prod_{\chi \in I} \BP(V_{\chi}(H))
$$ 
is Zariski locally trivial with fibre $\big[\prod_{\chi\in I} \G_{m,\chi}\big]/H$
this implies the second assertion in Claim \ref{claim2.1}.
\end{proof}
By \eqref{eq2.2} and \eqref{eq2.3}
$\displaystyle
[V/G]=\sum_{I\subset \hat{H}}\Big( [\G_m^{\#I}]\cdot \prod_{\chi\in I} [ \BP^{r_\chi}_k]\Big).$\\[.1cm]
This decomposition just depends on the dimensions $r_\chi+1$ of the subspaces
$V_{\chi}(H)$. So if $W_\chi$ denotes any $k$-vectorspace of this dimension and
$W=\bigoplus_{\chi\in \hat{H}}W_\chi$, one finds in $K_0({\rm Var}_k)$ 
$$
\BL^{\dim_K V}=\BL^{\dim_k W}= \sum_{I\subset \hat{H}} \ \prod_{\chi\in I} [W_\chi^\times]=
\sum_{I\subset \hat{H}}\Big( [\G_m^{\#I}]\cdot \prod_{\chi\in I} [ \BP^{r_\chi}_k]\Big)=
[V/G].
$$
This finishes the proof of Lemma \ref{thm1.1}.
\qed
\section{Verification of the properties in Example \ref{ex1.2}}
In the standard basis $e_1,e_2$ of $\Q^2$ and the basis $(1, \sqrt{-1})$ of $K/\Q$, we write 
\ga{}{\sigma: (x_1+\sqrt{-1} y_1)e_1+(x_2+ \sqrt{-1} y_2)e_2\mapsto (-x_1 +\sqrt{-1} y_1)e_2 +(x_2-\sqrt{-1} y_2)e_1,\notag}
As $\sigma$ is $\Q$-linear, 
it leaves the origin of $V$ invariant, thus acts on $V^\times=V\setminus \{0\}$. 
One has $\sigma^2=-{\rm Id}$ and this defines the extension
\ga{}{0\longrightarrow H:=\langle \sigma^2\rangle \longrightarrow G \longrightarrow \Gamma:=\langle \gamma \rangle \longrightarrow 0 \notag \\
\mbox{with \ \ } \Gamma=\langle \gamma \rangle\cong \Z/2={\rm Aut}(\Q(\sqrt{-1})/\Q), \mbox{ \ \ and \ \ } \gamma(\sqrt{-1})=-\sqrt{-1}. \notag}
Thus one has the fiber square
\ga{}{ \xymatrix{ \ar @{} [dr] |{\square}\ar[d] V/H \ar[r] & \Spec K \ar[d]\\
V/G \ar[r] & \Spec \Q} \notag
}
The  $\G_m$-bundle $V^\times \to \BP^1_K$ is compatible with the $G$-action. The subgroup $H$ acts trivially on $\BP^1_K$ while $\sigma$ acts via 
\ga{}{\bar \sigma: (x_1+\sqrt{-1} y_1:x_2+\sqrt{-1} y_2)\mapsto (x_2-\sqrt{-1} y_2:-x_1+\sqrt{-1}y_1). \notag}
This yields the fiber squares
\ga{}{
\xymatrix{ \ar @{} [dr] |{\square}\ar[d] V^\times /H \ar[r] & \ar @{} [dr] |{\square}
\ar[d] \ar[d]^{\pi} \BP^1_K \ar[r] & \Spec K \ar[d]\\
V^\times/G \ar[r] & \BP^1_K/G \ar[r] & \Spec \Q \notag}
}
\begin{claim} \label{claim3.1}
 $\BP^1_K/G$ is a genus 0 curve over $\Q$ without a rational point. 
\end{claim}
\begin{proof}
Indeed, a rational point is a fixpoint of $\BP^1_K$ under $\bar \sigma$. But the equation for a fixpoint is precisely 
$$
x^2_1+y^2_1+x_2^2+y_2^2=0, \mbox{ \ \ with \ \ } (x_1,x_2,y_1,y_2)\neq (0,0,0,0).$$ 
So over $\Q$ there are no solutions. 
\end{proof}
\begin{cor} \label{cor3.2}
 $\BL^2 \neq [V/G] \in K_0({\rm Var}_\Q)$. 
\end{cor}
\begin{proof}
The origin $x_1=x_2=y_1=y_2=0$ in $V$ is a fixpoint under $G$. 
Thus 
$$[V/G]=[V^\times/G] + [\Spec \Q].$$  
On the other hand, as we have seen in Claim \ref{claim2.1} $V^\times/G \to \BP^1_K/G$ is 
a locally trivial $\G_m$ bundle. \\[.1cm]
Here the trivialization of can be written down explicitly:
$V^\times $ is the total space of the $\G_m$-bundle to the invertible sheaf $\sO_{\BP^1_K}(-1)$, while $V^\times/H\to \BP^1_K$ is the total space of the $\G_m$-bundle to the invertible sheaf $\sO_{\BP^1_K}(-2)=\pi^*\sL$, where $\sL\in \Pic(\BP^1_K/G)$. So $V^\times/G\to \BP^1_K/G$ is the $\G_m$-bundle to the invertible sheaf $\sL$. \\[.1cm]
One concludes
\ga{}{ [V/G]-[\Spec \Q]= [V^\times/G]=[\G_m]\cdot [\BP^1_K/G]\in K_0({\rm Var}_{\Q}). \notag}
On the other hand, one also has
\ga{}{\BL^2-[\Spec \Q] = [\A^2_\Q\setminus \{0\}]=[\G_m]\cdot [\BP^1_\Q] \in K_0({\rm Var}_{\Q}). \notag}
If $[V/G]$ was equal to $\BL^2$ in $K_0({\rm Var}_\Q)$, then one would have the relation $[V^\times/G]=[\A^2_\Q\setminus \{0\}]$ in $K_0({\rm Var}_\Q)$, thus the relation
\ga{}{ \Phi([V^\times/G])=-\Phi([\BP^1_K/G])=\Phi([\A^2_\Q\setminus \{0\}])=-\Phi([\BP^1_\Q])  \mbox{ \ in \ }  \Z[SB],\notag}
where $\Phi: K_0({\rm Var}_{\Q}) \to \Z[SB]$ maps the class $[X]$ of a smooth projective $\Q$-variety $X$ to its stably birational equivalence class.\\[.1cm]
This however contradicts Claim \ref{claim3.1}, as the existence of a rational point is compatible with the stably birational equivalence on smooth projective varieties over any infinite field $k$.\\[.1cm]
For sake of completeness let us recall the proof of this well known fact. If $\tau: V\dashrightarrow W$ is a birational map between two smooth projective varieties, and $\tau $ is well defined near $v\in V(k)$, then $\tau(v)$ is well defined and lies in $W(k)$. Else one blows up $v$. This yields an exceptional divisor  $\BP^{{\rm dim}_k V-1}$. Since $\tau$ is well defined outside of codimension $\ge 2$,  and since $k$ is infinite, there are rational points on the exceptional divisor on which $\tau$ is defined and one repeats the argument. 
\end{proof}
\section{Proof of Proposition \ref{prop1.3}}
We first show 1). If $V$ has $k$-dimension $\le 2$,  we write the $G$-equivariant stratification $V=\{0\}\sqcup V^\times$. Furthermore, the projection $ V^\times\to \BP(V)$ is $G$-equivariant as well. Looijenga's argument shows here 
$$[V^\times/G]=[\G_m]\cdot [\BP(V)/G]\in K_0({\rm Var}_k).$$ 
On the other hand,  either 
$$\BP(V)=\Spec k=\BP(V)/G \mbox{ \ \ or \ \ }\BP(V)/G\cong_k \BP^1_k \cong_k \BP(V).$$ Adding up, one finds $[V/G]=\BL^2\in K_0({\rm Var}_k)$.\\[.2cm]
We now show 2). Instead of the decomposition $V=\bigoplus_{i=1}^r V_i\otimes T_i$ of $V$ as a direct sums of irreducible $G$ representations, considered in the introduction, we will
drop the condition that ${\rm Hom}_G(V_i,V_j)=\delta_{ij}\cdot k$ and choose
a decomposition $V=\bigoplus_{i=1}^m V_i$ as a direct sum of irreducible representations.
As usual we consider $V$ as a variety and write 
\eq{4.1}{V=\prod_{i=1}^m V_i.}
The monodromy group, that is the image of $G$ in $GL_k(V)$
is still a $p$-order cyclic group. So we may assume  
\eq{4.2}{G\subset GL_k(V)} in the discussion. 
\begin{claim} \label{claim4.1}
 There is a direct factor  $V_i$ of \eqref{4.1} such that $G\subset GL_k(V_i)$.
\end{claim}
\begin{proof}
 Since a $p$-power order cyclic group $G$ contains a unique $p$-order cyclic subgroup $C(G)$, if $\{1\}\neq K_i:={\rm Ker} \big(G\to GL_k(V_i)\big)$, then $C(G)= C(K_i)\subset K_i$. We conclude by \eqref{4.2}.
\end{proof}
We now change the notation: we set $U=V_i$ and $W=\bigoplus_{j\neq i} V_j$ with $V_i$ constructed in Claim \ref{claim4.1}. So $V=U\oplus W$
equivariantly. We assume that the dimension of $U$ is $2$. If this is $1$, the argument simplifies enormously and we don't detail. We define the $G$-equivariant stratifications
\ga{4.3}{U=\{0\}\sqcup D^\times \sqcup U^{(2)}\\
V=(\{0\}\times_k W)\sqcup (D^\times\times_k W)\sqcup (U^{(2)}\times_k W).\notag}
The strata are defined as follows. Write $\langle \sigma \rangle=G$. Let $F(T)\in k[T]$ be the minimal polynomial of $\sigma$. Since $U$ is irreducible, $F(T)$ is also the characteristic polynomial of $\sigma$ on $U$. This defines the quadratic extension 
\eq{4.4}{ K=k[T]/(F(T)).}
The linear map $\sigma\otimes K \in GL(U\otimes K)$ has two conjugate eigenlines and $$D=\{0\}\sqcup D^\times\subset U$$ is the $k$-irreducible curve defined by the union of the two lines. 
Further $$U^{(2)}=U\setminus D.$$ By definition, $G$ acts fixpoint free on $U^{(2)}$. 
\begin{claim} \label{claim4.2}
$ [(U^{(2)}\times_k W)/G]=[(U^{(2)}/G)\times_k W]=[U^{(2)}/G]\cdot [W]\in K_0({\rm Var}_k)$.
\end{claim}
\begin{proof}
 One has a $G$-equivariant projection $q: (U^{(2)}\times_k W)/G\to U^{(2)}/G$. Since $G\subset GL_k(U)$, for all points $x\in U^{(2)}$ with residue field $\kappa(x)\supset k$, one has $q^{-1}(x)\cong_{\kappa(x)} W\otimes_k \kappa(x)$. By construction, one has a fiber square
\eq{4.5}{
\xymatrix{ \ar @{} [dr] |{\square}\ar[d] U^{(2)}\times_k W \ar[r] & (U^{(2)}\times_k W)/G \ar[d]^q\\
U^{(2)} \ar[r] & U^{(2)}/G.}
}
Since $U^{(2)} \to  U^{(2)}/G$ is \'etale, $q$ defines a local system  in $H^1_{{\rm \acute{e}t}}(U^{(2)}/G, G_W)$ where $G_W$ is the image of $G$ in $GL_k(W)$. Then $(U^{(2)}\times_k W)/G$ is the total space of the  torsor in $H^1_{{\rm \acute{e}t}}(U^{(2)}/G, GL_k(W))$ induced by $G_W\inj GL_k(W)$. By flat descent \cite[Lemma~4.10]{Mi}, $$H^1_{\rm{\acute{e}t}}(U^{(2)}/G, GL_k(W))=H^1_{{\rm Zar}}(U^{(2)}/G, GL_k(W)).$$ 
Thus $(U^{(2)}\times_k W)/G\xrightarrow{q} U^{(2)}/G$, as the total space of a vector bundle, is Zariski locally trivial. We conclude
\ga{4.6}{ [(U^{(2)}\times_k W)/G]=[U^{(2)}/G]\cdot [W] \in K_0({\rm Var}_k).}
\end{proof}
So using \eqref{4.3} and Claim \ref{claim4.2}, we see
\ml{4.7}{ \ \ \ \ \ [V]-[V/G]= \big ([W]-[W/G]\big) +\\ \big([D^\times\times_k W]-[(D^\times\times_k W)/G] \big) + ([U^{(2)}]-[U^{(2)}/G])\cdot [W]. \ \ \ \ \ }
The curve $D^\times$ is $k$-irreducible, but splits over $K$. Therefore $K\subset H^0(D^\times, \sO)$ is the algebraic closure of $k$ and thus $G$ acts on $K$.
\begin{claim} \label{claim4.3}
 The action of $G$ on $\Spec K$ is trivial.
\end{claim}
\begin{proof}
After the choice of a cyclic vector, $\sigma$ is the matrix $\begin{pmatrix}
                                                              0 & 1\\b &a
                                                             \end{pmatrix}$
with $a,b\in k$. 
The curve $D^\times$ is $k$-affine. Its affine ring is $$H^0(D^\times, \sO)= k[X,Y, \frac{1}{X}]/\langle f(X,Y) \rangle $$ where the homogeneous polynomial $f(X,Y)=Y^2 -aXY -bX^2$ defines the irreducible polynomial $F(T)=T^2-aT-b$ yielding the $k$-quadratic extension   $K$. The inclusion of $K\subset H^0(D^\times, \sO)$
is $k$-linear and defined by $T\mapsto \frac{Y}{X}$. Furthermore, 
 $\sigma(X)=Y, \ \sigma(Y)=bX+aY$, thus 
$$
\sigma(T)= \frac{\sigma(Y)}{\sigma(X)}=\frac{bX+aY}{Y}= \frac{b}{T} +a=T.
$$
\end{proof}
We can now analyze the second difference in \eqref{4.7}.  
One has the $G$-equivariant  fiber product
\eq{}{
\xymatrix{ \ar @{} [dr] |{\square}\ar[d] D^\times\times_k W \ar[r] & \Spec K\times_k W \ar[d]\\
D^\times \ar[r] & \Spec K.} \notag
}
Since $D^\times =\Spec K\times_k \G_m$, the morphism $D^\times\times_k  W\to \Spec K\times_k W$ is a $G$-equivariant Zariski locally trivial $\G_m$-fibration. We first deduce
$$[D^\times\times_k  W]=[\G_m]\cdot [\Spec K]\cdot [W].$$
From the induced fiber square 
\eq{}{
\xymatrix{ \ar @{} [dr] |{\square}\ar[d] (D^\times\times_k W)/G \ar[r] & (\Spec K\times_k W)/G \ar[d]\\
(D^\times)/G \ar[r] & (\Spec K)/G=\Spec K} \notag
}
and $(D^\times)/G =(\Spec K\times_k \G_m)/G=\Spec K\times_k (\G_m/G)=\Spec K\times_k \G_m$,
 we deduce that 
$(D^\times\times_k  W)/G\to (\Spec K\times_k W)/G$ is a
Zariski locally trivial $\G_m$-fibration, and thus 
$$[(D^\times\times_k  W)/G]=[\G_m]\cdot [\Spec K]\cdot [W/G].$$ 
We conclude
\eq{4.8}{ [D^\times\times_k  W]-[(D^\times\times_k  W)/G]= [\G_m]\cdot [\Spec K]\cdot \big( [W]-[W/G]\big).}
We now analyze the third difference in \eqref{4.7}.
One has the  $G$-equivariant fibration $U^{(2)}\to \BP(U)\setminus \Spec K$, which is a $\G_m$-bundle. So 
$$[U^{(2)}]=[\G_m]\cdot ([\BP(U)]-[\Spec K]).$$ 
Since $\BP(U)/G$ is $k$-isomorphic to $\BP^1_k$, the group $G$ acts trivially on $\Spec K$, and  $U^{(2)}/G \to  \BP(U)/G $ is a $\G_m$-bundle,
one has 
\ml{4.9}{[U^{(2)}/G]=[\G_m]\cdot ([\BP(U)/G]-[\Spec K])=\\ [\G_m]\cdot ([\BP(U)]-[\Spec K])= 
[U^{(2)}]\in K_0({\rm Var}_k).}
Summing up, \eqref{4.7} reads 
\eq{4.10}{ [V]-[V/G]= \big(1+[\G_m]\cdot [\Spec K]\big)\cdot ([W]-[W/G]).}
Now $W$ has one less irreducible factor than $V$. We argue by induction on the number of irreducible factors, applying 1) to start the induction. This finishes the proof. \qed

\bibliographystyle{plain}

\begin{thebibliography}{99}
\bibitem{AKMW} Abramovich, D.; Karu, K.; Matsuki, K.; W\l odarczyk, 
J.: {\em Torification and factorization of birational maps},  J. Amer. Math. Soc. {\bf 15}  (2002), 531--572.
\bibitem{Bi}  Bittner, F.: {\em  The universal Euler characteristic for varieties of characteristic zero},  Compositio math., {\bf 140}(4) (2004), 1011--1032.
\bibitem{Ek} Ekedahl, T.:  {\em A geometric invariant of a finite group}, arXiv:0903.3148. 
\bibitem{Ko} Koll\'ar, J.: {\em  Conics in the Grothendieck ring},  Adv. Math., {\bf 198} (1) (2005), 27--35.
\bibitem{LaLu}  Larsen, M.,  Lunts, V. A.: {\em  Motivic measures and stable birational geometry}, Mosc. Math. J., {\bf 3}(1) (2003), 85--95.
\bibitem{Loo} Looijenga, E.: {\em Motivic measures}, S\'eminaire Bourbaki  1999/2000, Ast\'erisque {\bf 276} (2002), 267--297. 
\bibitem{Mi} Milne, J.; {\em \'Etale cohomology},  Princeton Mathematical Series {\bf 33} (1980), Princeton University Press. 
\bibitem{Ni} Nicaise, J.: {\em A trace formula for varieties over a discretely
valued field}, arXiv:0805.1323v2.
\bibitem{Poo} Poonen, B.: {\em The Grothendieck ring of varieties is not a domain}, Math. Res. Let. {\bf 9} (4) (2002), 493--498.
\bibitem{Ro} R\"okaeus, K.: {\em  The computation of the classes of some tori in the Grothendieck ring of varieties}, arXiv:0708.4396.
\bibitem{S} Saltman, D. J.: {\em Noether's problem over an algebraically closed field}, Inventiones mathematicae {\bf 77} (1984), no 1, 71--84.

\end{thebibliography}
\renewcommand\refname{References}

\end{document}